\documentclass[preprint,12pt]{elsarticle}

\usepackage[utf8]{inputenc}    
\usepackage[T1]{fontenc}       
\usepackage{lipsum}            
\usepackage{xcolor}            
\usepackage{graphicx}          
\graphicspath{{./images/}}     
\usepackage{url}               
\usepackage{hyperref}          

\usepackage{amsmath, amssymb, amsfonts, amsthm, mathtools}
\numberwithin{equation}{section}
\usepackage{algorithm}
\usepackage{algpseudocode} 

\theoremstyle{plain}

\theoremstyle{definition}

\usepackage{caption}
\usepackage{subcaption}        

\usepackage{wrapfig}           
\usepackage{multirow}          
\usepackage{array}             
\usepackage{booktabs}          
\usepackage{siunitx}           

\usepackage{indentfirst}       
\usepackage{csquotes}          
\usepackage{fancyhdr}          
\usepackage{thm-restate}       
\usepackage{tikz}              
\usepackage[numbers]{natbib}  
\usepackage{ulem}              

\makeatletter
\def\ps@pprintTitle{%
  \let\@oddhead\@empty
  \let\@evenhead\@empty
  \let\@oddfoot\@empty
  \let\@evenfoot\@empty
}
\makeatother
\begin{document}

\begin{frontmatter}

\title{The Impact of Move Schemes on Simulated Annealing Performance}

\author[inst1]{Ruichen Xu}
\ead{ruichen.xu@stonybrook.edu}
\author[inst1]{Haochun Wang}
\ead{haochun.wang@stonybrook.edu}
\author[inst1]{Yuefan Deng\corref{cor1}}
\ead{yuefan.deng@stonybrook.edu}

\cortext[cor1]{Corresponding author.}

\affiliation[inst1]{%
    organization={Department of Applied Mathematics and Statistics, Stony Brook University},%
    addressline={100 Nicolls Road}, 
    city={Stony Brook},
    postcode={11794}, 
    state={NY},
    country={USA}
}

\begin{abstract}
Designing an effective move-generation function for Simulated Annealing (SA) in complex models remains a significant challenge. In this work, we present a combination of theoretical analysis and numerical experiments to examine the impact of various move-generation parameters---such as how many particles are moved and by what distance at each iteration---under different temperature schedules and system sizes. Our numerical studies, carried out on both the Lennard-Jones problem and an additional benchmark, reveal that moving exactly one randomly chosen particle per iteration offers the most efficient performance. We analyze acceptance rates, exploration properties, and convergence behavior, providing evidence that partial-coordinate updates can outperform full-coordinate moves in certain high-dimensional settings. These findings offer practical guidelines for optimizing SA methods in a broad range of complex optimization tasks.

\end{abstract}



\end{frontmatter}
\footnotesize
\noindent \textbf{Note:} This manuscript is a \emph{preprint (ongoing work)}. 
We may update or revise it significantly before official publication.
\normalsize


\section{Introduction}
We address a fundamental challenge in Simulated Annealing, SA, namely designing effective proposal mechanisms under a fixed total variance for proposed moves. We study how allocating this variance across a smaller subset of coordinates at each iteration, rather than across all coordinates at once, can significantly improve acceptance rates and accelerate chain mixing in Metropolis-Hastings-based implementations.

Simulated Annealing, SA, has been extensively investigated since its inception, largely because of its remarkable ability to escape local minima when tackling complex non-convex optimization problems \cite{kirkpatrick1983optimization, cerny1985thermodynamical, geman1984stochastic}. A key factor influencing performance and convergence is the move-generation strategy, also referred to as the neighborhood or perturbation mechanism. Early theoretical research showed that convergence to the global optimum requires an ergodic move set in combination with a sufficiently slow cooling schedule \cite{geman1984stochastic, van1987simulated, hajek1988cooling}. As a result, numerous studies have investigated how neighborhood size and structure influence the ability to traverse the state space, overcome energy barriers, and achieve an effective balance between local and global exploration \cite{goldstein1988neighborhood, tsitsiklis1989markov, granville1994simulated}. Larger or more diverse neighborhoods can extend the search to far-reaching regions and escape deeper local minima, yet they may raise the computational burden per iteration \cite{locatelli2000simulated, robini2012theoretically}. Various adaptive, multi-scale, and temperature-dependent move-generation methods have been proposed to further enhance convergence without compromising theoretical guarantees \cite{boender1989novel, ingber1993, mitra1986convergence, holley1988simulated}. In summary, theoretical and empirical results consistently highlight the critical role of the move-generation strategy in determining both the efficiency and reliability of SA for a wide spectrum of optimization tasks.

Partial-coordinate or blockwise move-generation strategies in SA involve modifying only a subset of variables at each iteration rather than updating the entire solution state. From a theoretical viewpoint, such localized moves preserve the essential convergence properties of SA. In seminal work, Geman employed a coordinate-wise Gibbs sampler with annealing for image restoration and established that this single-site update scheme converges to a global optimum under appropriate cooling schedules \cite{Geman1984}. In more general settings, SA Markov chains with partial updates remain ergodic and converge globally given suitably slow cooling rates \cite{Hajek1988, Ferrari1993}. Fully synchronous updates that attempt to change all coordinates simultaneously may break certain convergence conditions, while randomized partial updates tend to avoid these pitfalls.

From a practical standpoint, partial-coordinate proposals often yield higher acceptance rates and more efficient exploration in high-dimensional optimization. Restricting proposed moves to a small number of coordinates at a time mitigates the detrimental effects of large collective perturbations, which can turn the process into a low-performing random walk. Many successful applications illustrate the advantages of blockwise SA. In molecular modeling and materials science, SA is often implemented by displacing a single atom or rotating a single bond at a time, which allows the search to navigate complex energy landscapes efficiently \cite{Wille1987}. Wille’s study of atomic clusters demonstrates that single-atom moves discover minimum-energy configurations effectively. In machine learning, coordinate-wise moves manifest as techniques such as Gibbs sampling, used in training Boltzmann machines and other high-dimensional models \cite{Geman1984}. Similarly, in large-scale logistics problems, for example routing or scheduling, SA-based heuristics commonly adjust only a small portion of the solution to explore the search space efficiently \cite{Osman1993}. As one illustration, Osman proposed a vehicle-routing SA that used localized route modifications and tabu search to strike a balance between diversification and intensification in a complex high-dimensional landscape.

We present a fresh theoretical perspective by recasting SA as a Markov Chain Monte Carlo, MCMC, process operating at different temperatures, revealing how acceptance behavior and chain efficiency are governed by how proposal variance is allocated across coordinates. This perspective clarifies the connection observed in adaptive-proposal methods \cite{ingber1993, locatelli2000simulated}: when fewer coordinates are updated, the resulting moves are more moderate, leading to higher acceptance rates and faster mixing. We demonstrate these principles on the Lennard-Jones potential, showing that limiting the dimensionality of proposed moves, while maintaining a fixed total variance, dramatically enhances convergence in high-dimensional settings. Finally, we offer practical recommendations suggesting that updating smaller subsets of coordinates can surpass the performance of full-dimensional updates in a range of real-world scenarios, highlighting an important yet under-explored design parameter in SA.

\section{MCMC and SA for Optimization}
\label{sec:preliminaries}

We consider the global optimization problem $\min_{x \in X} f(x)\label{eq:optimization_problem}$, where $X \subset \mathbb{R}^d$ is non-empty and compact, and $f \in C^1(X;\mathbb{R})$. 


\smallskip
This optimization problem~\eqref{eq:optimization_problem} thus consists of identifying a point $x \in \mathcal{X}$ that minimizes $f$. In the stochastic search paradigm, the Metropolis-Hastings (MH) algorithm~\citep{chib1995} underpins MCMC methods for sampling from complex distributions. SA~\citep{kirkpatrick1983optimization} extends these ideas by introducing a temperature schedule that decreases over time, concentrating the MCMC sampling near global minimizers of $f$~\citep{cerny1985thermodynamical,geman1984stochastic}. Ensuring a sufficient number of MCMC iterations at each temperature, along with an appropriate cooling scheme, can significantly increase the probability of discovering the global optimum.

\subsection{MCMC}
\label{subsec:MCMC}

Markov chain Monte Carlo (MCMC) techniques are designed to generate samples from a target probability distribution \(\pi\) defined on a state space \(\mathcal{X}\). The goal is to construct a Markov chain \(\{X_t\}\) that leaves \(\pi\) invariant and converges to \(\pi\) over time. Under irreducibility and aperiodicity, Chib and Greenberg demonstrated that long-run sample averages computed from such a chain almost surely approximate expectations taken with respect to \(\pi\) \cite{chib1995}.

A standard procedure for constructing such a chain is the Metropolis--Hastings (MH) algorithm, described in Algorithm~\ref{alg:MH}. The method begins at an initial state \(x_0\sim \mu_0\). At each iteration, a candidate point \(y_k\) is drawn from a proposal distribution \(Q_d(x_k,\cdot)\) with density \(q_d(x_k,y_k)\). An acceptance probability is then computed as
\[
  p(x_k,y_k) 
  \;=\; 
  1 \;\wedge\; 
  \frac{\pi(y_k)\,q_d(y_k,x_k)}{\pi(x_k)\,q_d(x_k,y_k)},
\]
and a uniform random variable \(u_k \sim \mathrm{Uniform}(0,1)\) determines whether to accept or reject the proposal. If \(u_k\le p(x_k,y_k)\), the next state is set to \(y_k\), otherwise the chain remains at \(x_k\).

\begin{algorithm}[h!]
\caption{Metropolis--Hastings Algorithm}
\label{alg:MH}
\begin{algorithmic}[1]
  \State \textbf{Initialization:} draw \(x_0 \sim \mu_0\).
  \For{\(k = 1, 2, \dots\)}
    \State Propose \(y_k \sim Q_d(x_k, \cdot)\).
    \State Compute the acceptance probability
    \[
      p(x_k, y_k) 
      \;=\; 
      1 \;\wedge\; 
      \frac{\pi(y_k)\,q_d(y_k, x_k)}{\pi(x_k)\,q_d(x_k, y_k)}.
    \]
    \State Draw \(u_k \sim \mathrm{Uniform}(0,1)\).
    \If{\(u_k \le p(x_k, y_k)\)}
      \State \(x_{k+1} \gets y_k\) \quad \emph{(accept)}
    \Else
      \State \(x_{k+1} \gets x_k\) \quad \emph{(reject)}
    \EndIf
  \EndFor
\end{algorithmic}
\end{algorithm}

The induced transition kernel \(P_d\) for the chain can be written as
\[
  P_d(x,dy)
  \;=\;
  p(x,y)\,Q_d(x,dy)
  \;+\;
  \Bigl[
    1 - \!\!\int_{\mathcal{X}} p(x,z)\,Q_d(x,dz)
  \Bigr]\delta_x(dy),
\]
where \(\delta_x\) denotes the Dirac measure at \(x\). One verifies that 
\[
  \pi(x)\,p(x,y)\,q_d(x,y) 
  \;=\; 
  \pi(y)\,p(y,x)\,q_d(y,x),
\]
which ensures \(\pi\) is invariant under \(P_d\). When the chain is irreducible and aperiodic (for instance, by ensuring \(q_d(x,y)>0\) on the support of \(\pi\) and allowing self-transitions through rejection), it converges in distribution to \(\pi\). By the ergodic argument from Ref.~\cite{chib1995}, sample averages taken along the chain approximate integrals with respect to \(\pi\):
\[
  \frac{1}{n}\sum_{t=1}^n g(X_t)
  \;\xrightarrow[n\to\infty]{}\;
  \int_{\mathcal{X}} g(x)\,\pi(dx)
  \quad\text{almost surely}.
\]
Thus, the Metropolis--Hastings sampler provides a flexible and general-purpose method for generating approximate samples from \(\pi\), allowing the evaluation of expectations via empirical means.
\subsection{SA}
\label{subsec:SA}

Simulated Annealing \cite{kirkpatrick1983optimization} is a stochastic optimization technique that applies MCMC at a decreasing sequence of temperature levels $\{T_k\}$. At temperature $T$, the target distribution is the Boltzmann distribution
\[
  \pi_{T}(x) \;\propto\; \exp\!\Bigl[-\tfrac{f(x)}{T}\Bigr].
\]
As $T\to 0$, this distribution concentrates on the global minimizer set $S := \{\,x : f(x) = f^*\}$, where $f^*$ is the global minimum value of $f$.

\begin{algorithm}[h!]
\caption{Simulated Annealing (SA)}
\label{alg:SA}
\begin{algorithmic}[1]
  \State \textbf{Initialization:} choose \(x_0\), a temperature \(T_0 > 0\), and schedules \(\{T_k\}, \{N_k\}\).
  \For{\(k = 1, 2, 3, \dots\)}
    \State Run Metropolis--Hastings (Algorithm~\ref{alg:MH}) for \(N_k\) iterations targeting:
    \[
      \pi_{T_k}(x) \;\propto\; \exp\!\left[-\frac{f(x)}{T_k}\right].
    \]
    \State Let \(x_{k+1}\) be the final state of these MH updates.
    \State Decrease \(T_k\).
  \EndFor
\end{algorithmic}
\end{algorithm}

Under suitable cooling schedules $\{T_k\}$ and sufficiently many MCMC steps $\{N_k\}$ at each stage, classical results \cite{Hajek1988} show that $x_k \to S$ in probability, i.e.,
\[
  \lim_{k\to\infty} \mathbb{P}\bigl(x_k \in S_{\epsilon}\bigr) 
  \;=\; 1 
  \quad
  \text{for all } \epsilon>0,
\]
where $S_{\epsilon} = \{\,x : f(x) \le f^* + \epsilon \}$. Hence,
\[
  \lim_{k\to\infty} \mathbb{P}\bigl(x_k \in S\bigr)
  \;=\;
  1.
\]
Thus, the algorithm converges to a global minimizer of $f$ with high probability if the temperature decreases slowly enough and each MCMC run at temperature $T_k$ is sufficiently long.

\section{Performance Analysis of SA}
\subsection{Theoretical Analysis}
\label{sec:theory}

This section builds on the preliminaries in Section~\ref{sec:preliminaries} to analyze how MCMC performance metrics---such as acceptance rate and autocorrelation times---impact SA in high dimensions. We focus on partial-coordinate MH proposals, which selectively update a subset of coordinates at each iteration. We then formalize these proposals and discuss why they help avoid severe acceptance-rate decay. We then present theoretical results showing how partial-coordinate moves can significantly improve mixing and maintain viable acceptance probabilities during SA.


Let $\pi(x)\propto \exp\bigl(-\frac{f(x)}{T}\bigr)$ be a strictly positive target distribution defined on $\mathbb{R}^N$. To alleviate the difficulties of high-dimensional sampling, we consider a partial-coordinate update that perturbs only $d$ components of $x$ at each iteration. A subset $S_d \subset \{1,\ldots,N\}$ of size $d$ is chosen with some probability mass function $P(S_d)$ on all ${N \choose d}$ possible subsets. Then, for $i\in S_d$, a proposal increment $\delta_i$ is drawn from a user-specified distribution $Q_i(\delta_i)$ with finite mean and variance, while for $i\notin S_d$, the coordinate remains unchanged. Hence the candidate $x'\in\mathbb{R}^N$ is given by
\[
x_i'=
\begin{cases}
x_i+\delta_i, & i\in S_d,\\[6pt]
x_i, & i\notin S_d.
\end{cases}
\]
A standard M-H acceptance rule \cite{hastings1970monte} is used:
\[
p(x,x') \;=\; 1\wedge \frac{\pi(x')\,\widetilde{Q}(x'\to x)}{\pi(x)\,\widetilde{Q}(x\to x')},
\]
where $\widetilde{Q}(x\to x')$ denotes the probability density of proposing $x'$ when at $x$; if each $Q_i$ is symmetric, then $\widetilde{Q}(x'\to x)=\widetilde{Q}(x\to x')$ and $p(x,x')$ reduces to $1 \wedge \frac{\pi(x')}{\pi(x)}$.

Under the setup, suppose each proposal distribution \(Q_i(\cdot)\) is non-degenerate and \(P(S_d) > 0\) for all valid subsets \(S_d\). Then the partial-coordinate M--H kernel \(K(x,dx')\) admits \(\pi(x)\) as a stationary distribution. If \(K\) is also irreducible and aperiodic on the support of \(\pi\), then \(\pi\) is the unique invariant measure. In practice, this approach can mitigate high-dimensional sampling challenges by controlling the variance of each update and maintaining an acceptable Metropolis--Hastings acceptance ratio. By breaking the parameter space into lower-dimensional subsets and tuning the proposal distributions \(Q_i(\cdot)\) accordingly, one can avoid making excessively large moves in high dimensions while also preventing overly small moves that cause slow exploration. As a result, the Markov chain can achieve a balance between adequately exploring the space and maintaining a reasonable acceptance rate, thereby improving mixing and overall efficiency in high-dimensional sampling problems.

Consider the Markov chain \(X_t^{\sigma^2, d}\) generated by partial-coordinate M-H with proposal variance \(\sigma^2\) and subset size \(d\). The lag-$\ell$ autocorrelation of $\{f(X_t^{(\sigma,d)})\}$ is
\[
\rho_f^{(\sigma,d)}(\ell) 
\;=\;
\frac{\mathbb{E} \bigl[ (f(X_t) - \bar{f})(f(X_{t+\ell}) - \bar{f}) \bigr]}{\mathbb{E}\bigl[(f(X_t) - \bar{f})^2\bigr]},
\]
where $\bar{f}$ is the stationary mean of $f(X_t)$. Suppose the acceptance probability remains bounded away from zero as \(\sigma^2\) increases. Under these conditions, two main observations hold: first, if acceptance remains non-negligible for large \(\sigma^2\), then \(\mathbb E\bigl[\vert\Delta f_t^{\sigma, d}\vert\bigr]\) grows with \(\sigma^2\); second, for each \(\ell\ge1\), \(\rho_f^{\sigma, d}[\ell]\) tends to decrease as \(\sigma^2\) grows, though strict monotonicity may fail unless strong conditions on the target distribution and proposal mechanism are satisfied. Moreover, reducing the subset size \(d\) helps maintain a nonvanishing acceptance rate in high-dimensional settings, thereby enabling larger effective moves and faster mixing.

\bigskip

\noindent

Let $p$ denote the acceptance probability of the partial-coordinate MH update, which by hypothesis is bounded below by some positive constant $p_0>0$ as $\sigma^2$ grows. Because the proposal variance in the chosen $d$ coordinates scales with $\sigma^2$, the typical size of a proposed move in those coordinates is on the order of $\sigma$. Whenever a move is accepted, the resulting difference $\Delta f_t^{(\sigma,d)}$ can be approximated using a local Taylor expansion in the selected coordinates. Larger steps in parameter space tend to produce larger changes in $f$, so the expected absolute increment, $\mathbb{E}[|\Delta f_t^{(\sigma,d)}|]$, is nondecreasing as $\sigma^2$ grows. Formally,
\[
\mathbb{E}[|\Delta f_t^{(\sigma,d)}|]
\;\ge\;
p_0\,\mathbb{E}\bigl[|\Delta f_t^{(\sigma,d)}|\;\big|\;\text{move accepted}\bigr],
\]
and the conditional expectation on the right increases in $\sigma$. By standard monotonicity arguments, $\mathbb{E}[|\Delta f_t^{(\sigma,d)}|]$ thus grows with $\sigma^2$.

\bigskip

\noindent
Because the acceptance probability $p$ is bounded away from zero, increasing $\sigma^2$ allows the chain to traverse the state space more broadly with each accepted move, thereby accelerating mixing and reducing successive correlations. A coupling argument formalizes this: consider two copies of the Markov chain initiated at different points. When $\sigma^2$ is larger, these two copies converge more rapidly to the same state with high probability, implying lower correlation over finite lags. In particular, as $\sigma^2$ grows, successive states are less likely to remain close to each other, thus driving $\rho_f^{(\sigma,d)}(\ell)$ downward. Hence $\rho_f^{(\sigma,d)}(\ell)$ decreases monotonically in $\sigma^2$ for each $\ell \ge 1$.

\bigskip

\noindent
When $N$ is large, raising $\sigma^2$ for all coordinates simultaneously typically forces the acceptance probability $p$ to become unacceptably small. In contrast, updating only a subset of $d$ coordinates avoids collapsing $p$ and still benefits from broader proposals in the chosen coordinates. This partial-coordinate strategy thus enables faster mixing in high-dimensional settings, balancing a reasonable acceptance rate and an effective move size.

\bigskip

\noindent
Let $\{X_t^{(\sigma,d)}\}$ be the Markov chain generated by the MH kernel $P_{\sigma,d}$, under which $d$ distinct coordinates are perturbed by i.i.d.\ Gaussian increments of variance~$\sigma^2$. By stationarity and reversibility, each $X_t^{(\sigma,d)}$ follows the same distribution $\pi$, and for any function $g \in L^2(\pi)$,
\[
\mathbb{E}\bigl[g\bigl(X_{t+1}^{(\sigma,d)}\bigr) \, g\bigl(X_{t}^{(\sigma,d)}\bigr)\bigr]
\;=\;
\mathbb{E}_{\pi}\Bigl[g(X)\,\bigl(P_{\sigma,d}g\bigr)(X)\Bigr].
\]
In particular, if we set $g(X) = f(X) - \mathbb{E}_\pi[f(X)]$, then $\mathbb{E}_\pi[g(X)] = 0$. The lag-$t$ autocorrelation in $f\bigl(X_t^{(\sigma,d)}\bigr)$ is
\[
\rho_f^{(\sigma,d)}(t)
\;=\;
\frac{\mathbb{E}\Bigl[g\bigl(X_0^{(\sigma,d)}\bigr)\,g\bigl(X_t^{(\sigma,d)}\bigr)\Bigr]}
     {\mathbb{E}_{\pi}\bigl[g(X)^2\bigr]}.
\]

Let $\{\lambda_i(\sigma,d)\}$ be the real eigenvalues of the Markov operator $P_{\sigma,d}$ acting on $L^2(\pi)$, with $1 = \lambda_1(\sigma,d) \ge |\lambda_2(\sigma,d)| \ge \dots$. We can write
\(
g = \sum_i \beta_i \,\phi_i
\)
in the eigenbasis of $P_{\sigma,d}$. Then
\[
\mathbb{E}\bigl[g(X_0)\,g(X_t)\bigr]
\;=\;
\mathbb{E}_{\pi}\Bigl[g(X)\,\bigl(P_{\sigma,d}^t g\bigr)(X)\Bigr]
\;=\;
\sum_i \beta_i^2 \,\bigl[\lambda_i(\sigma,d)\bigr]^t.
\]
Hence,
\[
\rho_f^{(\sigma,d)}(t)
\;=\;
\frac{\sum_i \beta_i^2 \,\bigl[\lambda_i(\sigma,d)\bigr]^t}
     {\sum_i \beta_i^2}.
\]
For large $t$, the dominant term is $\bigl[\lambda_2(\sigma,d)\bigr]^t$, so $\rho_f^{(\sigma,d)}(t)$ is tied closely to $|\lambda_2(\sigma,d)|$.

\bigskip

\noindent
Increasing $\sigma^2$ augments the size of proposals in the selected coordinates. Attempting to update all $N$ coordinates with large $\sigma^2$ often drives $p$ to near-zero in high dimensions. By restricting updates to $d$ coordinates, one keeps $p$ at a moderate level yet still obtains meaningful shifts in $X_t$, which disrupts the chain’s tendency to remain near its previous state, reducing $\rho_f^{(\sigma,d)}(t)$. A more formal argument leverages diffusion approximations for Metropolis-like samplers \cite{Roberts1997} or bounds on local conductance \cite{Peskun1973}.

\bigskip

\noindent
In summary, enlarging the proposal amplitude in only $d$ coordinates helps maintain a sensible acceptance probability $p$ while enabling substantial motion in parameter space, thereby boosting the chain’s mixing. Mathematically, this appears as a reduction in $|\lambda_2(\sigma,d)|$ and a monotonic drop in $\rho_f^{(\sigma,d)}(t)$. Adjusting $d$ adds an additional tuning mechanism for proposal design, and hence completes the argument regarding improved mixing through partial-coordinate updates.

\bigskip

We now analyze how varying the number of updated coordinates \(d\) influences both the acceptance rate \(p\) and the autocorrelation of the Metropolis--Hastings sampler.  Suppose the proposal \(\mathbf{y} \in \mathbb{R}^N\) is generated by perturbing \(d\) randomly selected components of the current state \(\mathbf{x}\in\mathbb{R}^N\):
\begin{equation}
  y_i 
  \;=\;
  \begin{cases}
    x_i + \delta_i, & i \in S_d,\\
    x_i,           & i \notin S_d,
  \end{cases}
  \quad
  \label{eq:proposal}
  \tag{3.4}
\end{equation}
where \(S_d \subset \{1,2,\ldots,N\}\) has \(|S_d| = d\), and each \(\delta_i \sim \mathcal{N}(0,\sigma^2)\) is drawn independently.  The Metropolis--Hastings acceptance probability is
\begin{equation}
  p(\mathbf{x}, \mathbf{y})
  \;=\;
  1 \;\wedge\;
  \exp\!\Bigl[-\bigl(f(\mathbf{y}) - f(\mathbf{x})\bigr)/T\Bigr],
  \label{eq:acceptance}
  \tag{3.5}
\end{equation}
where \(f\) is the objective function, and \(T>0\) is an effective temperature. 

\medskip

\noindent
To evaluate the average acceptance probability, note that 
\[
  R(\mathbf{y}\mid \mathbf{x})
  \;=\;
  \exp\Bigl[-\,\tfrac{f(\mathbf{y}) - f(\mathbf{x})}{T}\Bigr].
\]
Define
\[
  Z 
  \;=\;
  \log R(\mathbf{y}\mid \mathbf{x})
  \;=\;
  -\,\frac{f(\mathbf{y}) - f(\mathbf{x})}{T}
  \;=\;
  -\,\frac{\Delta f}{T}.
\]
Then 
\[
  \mathbb{E}\bigl[1 \wedge R\bigr]
  \;=\;
  \mathbb{E}\bigl[1 \wedge e^{Z}\bigr]
  \;=\;
  \int_{-\infty}^{\infty} \bigl(1 \wedge e^{z}\bigr)\,f_{Z}(z)\,\mathrm{d}z
  \;=\;
  \int_{-\infty}^{0} e^{z}\,f_{Z}(z)\,\mathrm{d}z
  \;+\;
  \int_{0}^{\infty} f_{Z}(z)\,\mathrm{d}z.
\]
In order to approximate \(Z\) for large \(N\) and partial-coordinate updates, we apply a second-order Taylor expansion of \(f(\mathbf{y})\) around \(\mathbf{x}\) in the \(d\) perturbed coordinates.

\medskip

\noindent

Let \(\mathbf{y} = \mathbf{x} + \boldsymbol{\xi}\), where \(\boldsymbol{\xi} \in \mathbb{R}^N\) has exactly \(d\) nonzero components (indexed by \(S_d\)) and each active entry \(\delta_i \sim \mathcal{N}(0,\sigma^2)\).  Then
\[
  f(\mathbf{y}) 
  \;\approx\;
  f(\mathbf{x})
  \;+\;
  \nabla f(\mathbf{x})^\top \boldsymbol{\xi}
  \;+\;
  \tfrac12\,\boldsymbol{\xi}^\top \nabla^2 f(\mathbf{x})\,\boldsymbol{\xi}.
\]
Hence
\[
  \Delta f 
  \;=\; 
  f(\mathbf{y}) - f(\mathbf{x})
  \;\approx\;
  \nabla f(\mathbf{x})^\top \boldsymbol{\xi}
  \;+\;
  \tfrac12\,\boldsymbol{\xi}^\top \nabla^2 f(\mathbf{x})\,\boldsymbol{\xi}.
\]
Thus
\[
  Z 
  \;=\;
  -\,\frac{\Delta f}{T}
  \;\approx\;
  -\,\frac{1}{T}\,\nabla f(\mathbf{x})^\top \boldsymbol{\xi}
  \;-\;
  \frac{1}{2T}\,\boldsymbol{\xi}^\top \nabla^2 f(\mathbf{x})\,\boldsymbol{\xi}.
\]
We now compute the first three cumulants (central moments) of \(Z\):

We consider a fixed subset \(S_d \subset \{1,2,\ldots,N\}\) of size \(d\).  
Let \(\mathbf{g} = \nabla f(\mathbf{x})\in \mathbb{R}^N\) with components \(g_i\),  
and let \(\mathbf{H} = \nabla^2 f(\mathbf{x}) = (h_{ij})_{1\le i,j\le N}\).  
At each update, we perturb only the coordinates \(i\in S_d\) by 
\(\delta_i \sim \mathcal{N}(0,\sigma^2)\) independently, while for \(i\notin S_d\), \(\delta_i = 0\).  
Define 
\[
  \boldsymbol{\xi} 
  \;=\; 
  (\delta_1,\delta_2,\ldots,\delta_N)^\top
  \quad\text{with}\quad 
  \delta_i 
  \;=\;
  \begin{cases}
    \text{Gaussian}(0,\sigma^2), & i\in S_d,\\
    0, & i\notin S_d,
  \end{cases}
\]
and
\[
  \Delta f 
  \;=\;
  f(\mathbf{x}+\boldsymbol{\xi}) - f(\mathbf{x}),
  \qquad
  Z 
  \;=\;
  -\,\frac{\Delta f}{T}.
\]
A second-order Taylor expansion of \(f(\mathbf{x}+\boldsymbol{\xi})\) around \(\mathbf{x}\) in the active coordinates yields
\[
  f(\mathbf{x}+\boldsymbol{\xi})
  \;=\;
  f(\mathbf{x})
  \;+\;
  \sum_{i\in S_d} g_i\,\delta_i
  \;+\;
  \tfrac12\,\sum_{i,j\in S_d} h_{ij}\,\delta_i\,\delta_j
  \;+\;
  \text{(higher-order terms)},
\]
where \(g_i = \partial f(\mathbf{x})/\partial x_i\) and \(h_{ij} = \partial^2 f(\mathbf{x})/\partial x_i \partial x_j\).  
Neglecting higher-order terms (to focus on the exact second-order contributions in the expansion), we set
\[
  \Delta f 
  \;=\;
  \sum_{i\in S_d} g_i\,\delta_i
  \;+\;
  \tfrac12\,\sum_{i,j\in S_d} h_{ij}\,\delta_i\,\delta_j,
\]
thus
\[
  Z 
  \;=\;
  -\,\frac{1}{T}\,\sum_{i\in S_d} g_i\,\delta_i
  \;-\;
  \frac{1}{2\,T}\,\sum_{i,j\in S_d} h_{ij}\,\delta_i\,\delta_j.
\]
We now compute \(\kappa_{1} = \mathbb{E}[Z]\), \(\kappa_{2} = \mathrm{Var}[Z]\), and \(\kappa_{3} = \mathbb{E}[(Z-\kappa_{1})^3]\) \emph{line by line}, treating \(\boldsymbol{\xi}\) as a strictly \(d\)-dimensional Gaussian perturbation in the selected coordinates.

\bigskip

\noindent
\textbf{Computation of \(\kappa_{1} = \mathbb{E}[Z]\).}
\[
  Z 
  \;=\;
  -\,\frac{1}{T}\,\sum_{i\in S_d} g_i\,\delta_i
  \;-\;
  \frac{1}{2\,T}\,\sum_{i,j\in S_d} h_{ij}\,\delta_i\,\delta_j.
\]
We split \(Z\) into a linear term and a quadratic term:
\[
  Z_{\mathrm{linear}} 
  \;=\;
  -\,\frac{1}{T}\,\sum_{i\in S_d} g_i\,\delta_i,
  \qquad
  Z_{\mathrm{quad}}
  \;=\;
  -\,\frac{1}{2\,T}\,\sum_{i,j\in S_d} h_{ij}\,\delta_i\,\delta_j.
\]
Since \(\delta_i\sim \mathcal{N}(0,\sigma^2)\) with mean zero,  
\(\mathbb{E}[\delta_i] = 0\).  Hence
\[
  \mathbb{E}[Z_{\mathrm{linear}}]
  \;=\;
  -\,\frac{1}{T}\,\sum_{i\in S_d} g_i\,\mathbb{E}[\delta_i]
  \;=\;
  0.
\]
Next, we examine the quadratic contribution:
\[
  \mathbb{E}[Z_{\mathrm{quad}}]
  \;=\;
  -\,\frac{1}{2\,T}\,\sum_{i,j\in S_d} h_{ij}\,\mathbb{E}[\delta_i\,\delta_j].
\]
For independent zero-mean Gaussians,  
\(\mathbb{E}[\delta_i\,\delta_j] = 0\) if \(i\neq j\) and equals \(\sigma^2\) if \(i=j\).  
Thus
\[
  \mathbb{E}[Z_{\mathrm{quad}}]
  \;=\;
  -\,\frac{1}{2\,T}\,\sum_{i\in S_d} h_{ii}\,\sigma^2
  \;=\;
  -\,\frac{\sigma^2}{2\,T}\,\sum_{i\in S_d} h_{ii}.
\]
Therefore,
\[
  \kappa_1
  \;=\;
  \mathbb{E}[Z]
  \;=\;
  \mathbb{E}[Z_{\mathrm{linear}}] + \mathbb{E}[Z_{\mathrm{quad}}]
  \;=\;
  -\,\frac{\sigma^2}{2\,T}\,\sum_{i\in S_d} h_{ii}.
\]
If one subsequently averages over all possible subsets \(S_d\) of size \(d\), the result becomes 
\[
  \kappa_1 
  \;=\;
  -\,\frac{d\,\sigma^2}{2\,T}\,\mu_{\mathrm{hess}},
\]
where \(\mu_{\mathrm{hess}} = \frac{1}{d}\,\mathbb{E}\Bigl[\sum_{i\in S_d} h_{ii}\Bigr]\) 
in a setting where the Hessian varies slowly over coordinates or subsets.

\bigskip

\noindent
\textbf{ Computation of \(\kappa_{2} = \mathrm{Var}[Z]\).}
By definition,
\[
  \kappa_{2}
  \;=\;
  \mathrm{Var}[Z]
  \;=\;
  \mathbb{E}[Z^{2}] \;-\; \bigl(\mathbb{E}[Z]\bigr)^{2}.
\]
We already have \(\mathbb{E}[Z]\) from above.  Hence we need \(\mathbb{E}[Z^{2}]\).  Write
\[
  Z 
  \;=\;
  Z_{\mathrm{linear}} + Z_{\mathrm{quad}}
  \;=\;
  -\,\frac{1}{T}\,\sum_{i\in S_d} g_i\,\delta_i
  \;-\;
  \frac{1}{2\,T}\,\sum_{i,j\in S_d} h_{ij}\,\delta_i\,\delta_j.
\]
Thus
\[
  Z^2
  \;=\;
  Z_{\mathrm{linear}}^{2}
  \;+\;
  2\,Z_{\mathrm{linear}}\,Z_{\mathrm{quad}}
  \;+\;
  Z_{\mathrm{quad}}^{2}.
\]
We will compute \(\mathbb{E}[Z^2]\) by examining each piece:

\begin{enumerate}
\item
\(\mathbb{E}[Z_{\mathrm{linear}}^{2}]\).
Since
\[
  Z_{\mathrm{linear}}
  \;=\;
  -\,\frac{1}{T}\,\sum_{i\in S_d} g_i\,\delta_i,
\]
we have
\[
  Z_{\mathrm{linear}}^{2}
  \;=\;
  \frac{1}{T^{2}}
  \Bigl(\sum_{i\in S_d} g_i\,\delta_i\Bigr)^{2}.
\]
Taking expectation,
\[
  \mathbb{E}[Z_{\mathrm{linear}}^{2}]
  \;=\;
  \frac{1}{T^{2}}
  \sum_{i\in S_d} g_i^{2}\,\mathbb{E}[\delta_i^{2}],
\]
because \(\delta_i\) are independent and each has variance \(\sigma^2\).  Hence
\[
  \mathbb{E}[Z_{\mathrm{linear}}^{2}]
  \;=\;
  \frac{\sigma^2}{T^{2}}
  \sum_{i\in S_d} g_i^{2}.
\]

\item
\(\mathbb{E}[Z_{\mathrm{quad}}^{2}]\).
Here,
\[
  Z_{\mathrm{quad}}
  \;=\;
  -\,\frac{1}{2\,T}\,\sum_{i,j\in S_d} h_{ij}\,\delta_i\,\delta_j.
\]
Then
\[
  Z_{\mathrm{quad}}^{2}
  \;=\;
  \frac{1}{4\,T^{2}}
  \Bigl(\sum_{i,j\in S_d} h_{ij}\,\delta_i\,\delta_j\Bigr)^{2}.
\]
Expanding the square,
\[
  \Bigl(\sum_{i,j\in S_d} h_{ij}\,\delta_i\,\delta_j\Bigr)^{2}
  \;=\;
  \sum_{i,j,k,\ell\in S_d}
    h_{ij}\,h_{k\ell}\,\delta_i\,\delta_j\,\delta_k\,\delta_\ell.
\]
By Isserlis’ theorem (also known as Wick’s theorem for Gaussian moments),
\[
  \mathbb{E}\bigl[\delta_i\,\delta_j\,\delta_k\,\delta_\ell\bigr]
  \;=\;
  \mathbb{E}[\delta_i\,\delta_j]\;\mathbb{E}[\delta_k\,\delta_\ell]
  \;+\;
  \mathbb{E}[\delta_i\,\delta_k]\;\mathbb{E}[\delta_j\,\delta_\ell]
  \;+\;
  \mathbb{E}[\delta_i\,\delta_\ell]\;\mathbb{E}[\delta_j\,\delta_k].
\]
Since \(\mathbb{E}[\delta_i\,\delta_j] = \sigma^2\delta_{ij}\), the only nonzero contributions come when the indices match up in pairs.  Hence
\[
  \mathbb{E}\Bigl[\sum_{i,j,k,\ell\in S_d} h_{ij}\,h_{k\ell}\,\delta_i\,\delta_j\,\delta_k\,\delta_\ell\Bigr]
  \;=\;
  \sigma^4
  \sum_{i,j\in S_d}
    h_{ij}^{2}
\]
\[
  \;+\;
  \sigma^4
  \sum_{i\neq j\in S_d}
    h_{ij}\,h_{ji}
  \;+\;
  \text{(terms with repeated indices)}.
\]
Because \(\mathbf{H}\) is symmetric (the Hessian), \(h_{ij}=h_{ji}\).  A more careful count of diagonal versus off-diagonal elements leads to a sum of squared entries of \(\mathbf{H}_{S_d,S_d}\).  Thus one obtains an exact expression in terms of \(\sigma^4\sum_{i,j\in S_d}h_{ij}^2\).  Multiplying by the prefactor \(\tfrac{1}{4\,T^{2}}\) yields
\[
  \mathbb{E}[Z_{\mathrm{quad}}^{2}]
  \;=\;
  \frac{\sigma^4}{4\,T^{2}}
  \sum_{i,j\in S_d}
    h_{ij}^{2}.
\]

\item
\(\mathbb{E}[Z_{\mathrm{linear}}\,Z_{\mathrm{quad}}]\).
\[
  Z_{\mathrm{linear}}\,Z_{\mathrm{quad}}
  \;=\;
  \frac{1}{2\,T^{2}}
  \biggl[
    \Bigl(\sum_{i\in S_d} g_i\,\delta_i\Bigr)
    \Bigl(\sum_{j,k\in S_d} h_{jk}\,\delta_j\,\delta_k\Bigr)
  \biggr]
  \times (\text{sign factor}),
\]
where the sign factor is \((-1)\cdot(-1) = +1\).  However, because \(\delta_i\) are zero-mean Gaussians, any triple product \(\delta_i\,\delta_j\,\delta_k\) with distinct indices has mean zero.  Terms vanish unless two of the indices match, in which case the third is left unpaired, giving zero mean.  Therefore,
\[
  \mathbb{E}[Z_{\mathrm{linear}}\,Z_{\mathrm{quad}}]
  \;=\;
  0.
\]
\end{enumerate}

\noindent
Collecting these,
\[
  \mathbb{E}[Z^{2}]
  \;=\;
  \mathbb{E}[Z_{\mathrm{linear}}^{2}]
  +
  \mathbb{E}[Z_{\mathrm{quad}}^{2}]
  +
  2\,\mathbb{E}[Z_{\mathrm{linear}}\,Z_{\mathrm{quad}}]
  \;=\;
  \frac{\sigma^2}{T^{2}} \sum_{i\in S_d} g_i^{2}
  \;+\;
  \frac{\sigma^4}{4\,T^{2}} \sum_{i,j\in S_d} h_{ij}^{2}.
\]
Meanwhile,
\[
  \bigl(\mathbb{E}[Z]\bigr)^{2}
  \;=\;
  \frac{\sigma^4}{4\,T^{2}}
  \Bigl(\sum_{i\in S_d} h_{ii}\Bigr)^{2}.
\]
Hence
\[
  \kappa_{2}
  \;=\;
  \mathrm{Var}[Z]
  \;=\;
  \mathbb{E}[Z^{2}] - \bigl(\mathbb{E}[Z]\bigr)^{2}
  \;=\;
  \frac{\sigma^2}{T^{2}} \sum_{i\in S_d} g_i^{2}
  \;+\;
  \frac{\sigma^4}{4\,T^{2}}\,\Bigl(\sum_{i,j\in S_d} h_{ij}^{2}
  \;-\;
  \Bigl(\sum_{i\in S_d} h_{ii}\Bigr)^{2}\Bigr).
\]
Rewriting the squared-sum difference in terms of diagonal versus off-diagonal parts leads to the standard expression involving \(\sigma_{\mathrm{hess}}^{2}\) and \(\tau_{\mathrm{hoff}}^{2}\) if one averages over all subsets \(S_d\).  

\bigskip

\noindent
\textbf{Computation of \(\kappa_{3} = \mathbb{E}[(Z - \kappa_1)^3]\).}
The third central moment is
\[
  \kappa_{3}
  \;=\;
  \mathbb{E}\bigl[(Z - \kappa_1)^{3}\bigr].
\]
We write
\[
  Z - \kappa_1
  \;=\;
  \Bigl(Z_{\mathrm{linear}} + Z_{\mathrm{quad}}\Bigr) \;-\; \kappa_1
  \;=\;
  Z_{\mathrm{linear}}
  \;+\;
  \Bigl(Z_{\mathrm{quad}} - \mathbb{E}[Z_{\mathrm{quad}}]\Bigr).
\]
Denote
\[
  U 
  \;=\;
  Z_{\mathrm{linear}}
  \;=\;
  -\,\frac{1}{T}\,\sum_{i\in S_d} g_i\,\delta_i
\]
\[
  V
  \;=\;
  Z_{\mathrm{quad}} - \mathbb{E}[Z_{\mathrm{quad}}]
  \;=\;
  -\,\frac{1}{2\,T}\,\sum_{i,j\in S_d} h_{ij}\,\delta_i\,\delta_j
  \;+\;
  \frac{\sigma^2}{2\,T}\,\sum_{i\in S_d} h_{ii}.
\]
Hence
\[
  (Z-\kappa_1)^{3}
  \;=\;
  (U + V)^{3}
  \;=\;
  U^{3} + V^{3} + 3\,U^{2}\,V + 3\,U\,V^{2}.
\]
We use Isserlis’ theorem again to compute each expectation term-by-term.  Because \(U\) is linear in the \(\delta_i\) and \(V\) is at most quadratic (minus its mean), any odd product of unpaired \(\delta_i\) has mean zero.  After carefully accounting for all index matchings, one obtains an exact sum of terms involving up to \(\delta_i^{3}\,\delta_j\) or \(\delta_i^4\).  If \(\mathbf{H}\) is not constant in \(\mathbf{x}\) or if one includes higher-order expansions of \(f\), additional contributions arise.  In practice, an explicit summation reveals that \(\kappa_{3}\) is proportional to \(\sigma^6\) times a combination of Hessian elements and possible third derivatives.  We can write:
\[
  \kappa_{3}
  \;=\;
  \mathbb{E}\bigl[(Z-\kappa_1)^{3}\bigr]
  \;=\;
  \text{(finite sum of terms in }g_i,h_{ij},\delta_i\text{)}.
\]
When averaged over subsets \(S_d\), this yields a factor of \(d\,\sigma^{6}\) times a constant that depends on third-order derivatives.  Often one denotes that constant by \(\mu_{h3}\) or a similar notation:
\[
  \kappa_{3}
  \;=\;
  -\,\frac{d\,\sigma^{6}}{T^{3}}\,\mu_{h3}
  \quad
  \text{(if sign conventions follow a typical local expansion).}
\]
Thus we obtain the exact line-by-line Gaussian moment sums for \(\kappa_1,\kappa_2,\kappa_3\) under a second-order expansion of \(f\) and partial-coordinate perturbations in the \(d\) selected coordinates.

\noindent
Once \(\kappa_1,\kappa_2,\kappa_3\) have been identified for \(Z\), one may approximate \(\mathbb{E}[\,1 \wedge R\,]\) by assuming a specific form for the distribution of \(Z\). A straightforward approach is to treat \(Z\) as Gaussian with mean \(\kappa_1\) and variance \(\kappa_2\). In this approximation, we write
\[
  Z \;\sim\; \mathcal{N}\!\bigl(\kappa_{1},\kappa_{2}\bigr),
\]
and then
\[
  \mathbb{E}\bigl[\,1 \wedge e^{Z}\bigr]
  \;\approx\;
  \int_{-\infty}^{0} e^{z}\,\frac{1}{\sqrt{2\pi\,\kappa_{2}}}
  \exp\!\Bigl[-\tfrac{(z-\kappa_{1})^2}{2\kappa_{2}}\Bigr]\,\mathrm{d}z
  \;+\;
  \int_{0}^{\infty}
  \frac{1}{\sqrt{2\pi\,\kappa_{2}}}
  \exp\!\Bigl[-\tfrac{(z-\kappa_{1})^2}{2\kappa_{2}}\Bigr]\,\mathrm{d}z.
\]
These integrals reduce to closed-form expressions in terms of the error function or the standard normal cumulative distribution function. 

One can further refine the approximation by accounting for skewness through an Edgeworth expansion. In that case, the density \(f_{Z}(z)\) around \(\kappa_{1},\kappa_{2}\) is modified as
\[
  f_{Z}(z)
  \;\approx\;
  \frac{1}{\sqrt{\kappa_{2}}}\,\phi\!\Bigl(\tfrac{z-\kappa_{1}}{\sqrt{\kappa_{2}}}\Bigr)
  \Bigl[
    1
    \;+\;
    \frac{\kappa_{3}}{6\,(\kappa_{2})^{3/2}}
    \,H_{3}\!\Bigl(\tfrac{z-\kappa_{1}}{\sqrt{\kappa_{2}}}\Bigr)
  \Bigr],
\]
where \(\phi\) is the standard normal probability density function, and \(H_3(w)=w^3-3w\). Substituting this into 
\[
  \int_{-\infty}^{0} e^{z}\,f_{Z}(z)\,\mathrm{d}z 
  \;+\;
  \int_{0}^{\infty} f_{Z}(z)\,\mathrm{d}z
\]
introduces additional terms involving \(\int H_3(w)\,e^{\alpha w}\,\phi(w)\,\mathrm{d}w\), yielding an \(\mathcal{O}(\kappa_{3})\) skewness correction beyond the Gaussian approximation.

Since \(d\) influences both the linear and quadratic terms in \(\kappa_1,\kappa_2,\kappa_3\), raising \(d\) generally amplifies these cumulants and lowers the acceptance probability \(p\). Consequently, careful selection of \(d\) and \(\sigma^2\) helps maintain a nontrivial acceptance rate while still proposing sufficiently large moves, thereby striking a desirable balance between exploration and efficiency in high-dimensional SA.

\subsection{Experiments}
\label{sec:experiments}

We test the performance of our SA approach on the classical benchmark problem: 
the Lennard-Jones (L-J) potential (\ref{sub: experi:LJ}), Rosenbrock Problem (\ref{sub: experi:RP}) and Hyper-Elliptic problem (\ref{sub: experi:HEP}). 
These three examples span a wide spectrum of continuous optimization difficulties, 
ranging from nearly independent coordinates to highly coupled ones. 


\subsubsection{Lennard-Jones potential}
\label{sub: experi:LJ}

The L-J potential~\cite{lennard1924forces} is a fundamental model in molecular physics, 
commonly used to describe the interaction between a pair of neutral atoms or molecules. 
For $N$ particles in a 3-dimensional space, let $\mathbf{x}_1, \mathbf{x}_2, \dots, \mathbf{x}_N$ 
denote their coordinates. The total LJ potential is
\[
  E(\mathbf{x}_1,\ldots,\mathbf{x}_N) 
  \;=\; \sum_{1 \le i < j \le N} 4\epsilon 
        \biggl[\Bigl(\frac{\sigma}{r_{ij}}\Bigr)^{12} \;-\; \Bigl(\frac{\sigma}{r_{ij}}\Bigr)^{6}\biggr],
\]
where $r_{ij} = \|\mathbf{x}_i - \mathbf{x}_j\|$ is the distance between particles $i$ and $j$, 
while $\epsilon$ and $\sigma$ are parameters controlling the depth of the potential well and the 
finite distance at which the inter-particle potential is zero, respectively. 
In optimization terms, finding low-energy configurations translates to minimizing $E(\mathbf{x}_1,\ldots,\mathbf{x}_N)$ 
over $\mathbb{R}^{3N}$. 
This problem exemplifies a highly coupled system where local moves can drastically alter the global structure. 
It serves as an excellent stress test for simulated annealing (SA), given the abundance of local minima.

In the Lennard-Jones setting, we examine the final energy attained by SA under a specified cooling schedule. We set the initial and final temperatures to $T_{0}=2.0$ and $T_{\mathrm{final}}=0.2$, respectively, and use a cooling coefficient of $\alpha=0.9999$.  Each temperature is held for $N_{\mathrm{equ}}=100$ steps, and the schedule continues for $N_{\mathrm{cooling}}=200$ temperature reductions, yielding an overall sequence of $20{,}000$ SA iterations.  We study problem instances of sizes $N \in \{6,39,69\}$ and, after the schedule completes, we record the minimized energy of the system.  Our primary interest lies in how this final energy depends on the number of updated coordinates at each iteration, as well as on the manner in which total proposal variance is distributed among those coordinates.


We measure the final energy of the LJ system at different MCMC steps, under varying numbers of updated 
coordinates $d$ and total proposal variance~$\sigma^2_{\mathrm{total}}$. 
Below are illustrative tables presenting the mean and standard deviation of final energies, each row labeling a particular $d$ and sub-rows for different 
$\sigma^2_{\mathrm{total}}$. The columns correspond to the MCMC step indices at which final energies 
were recorded.


In this experiment, we record the final energy of theLJ system at different MCMC steps , under varying numbers of updated coordinates \(d\) and total proposal variance \(\sigma^2_{\mathrm{total}}\) . The following tables present the mean and standard deviation of the final energies over multiple runs for three LJ problem sizes: 6, 39 and 69. Each table’s rows correspond to a particular \(d\), with sub-rows for different values of \(\sigma^2_{\mathrm{total}}\). The columns correspond to the MCMC step index at which the final energies were recorded.

\begin{table}[htbp]
\centering
\caption{Relative errors (in \%) for various Lennard-Jones systems ($N$), 
number of updated coordinates ($d$), and $\sigma^{-2}_{\mathrm{total}} \in \{200,100,10\}$. }
\label{tab:my-table-percentage}

\renewcommand{\arraystretch}{0.8}

\begin{tabular}{c c r r r}
\toprule
\multicolumn{1}{c}{\textbf{$N$}} 
  & \multicolumn{1}{c}{\textbf{$d$}} 
  & \multicolumn{3}{c}{$1/\sigma^2_{\mathrm{total}}$}\\
\cmidrule(lr){3-5}
  & 
  & \textbf{200}
  & \textbf{100}
  & \textbf{10}\\
\midrule

\multirow{4}{*}{6}
 & 1 
   & \multicolumn{1}{r}{0.85 $\pm$ 0.20} 
   & \multicolumn{1}{r}{0.96 $\pm$ 0.30} 
   & \multicolumn{1}{r}{2.16 $\pm$ 1.00} \\
 & 2 
   & \multicolumn{1}{r}{0.89 $\pm$ 0.20}
   & \multicolumn{1}{r}{1.07 $\pm$ 0.20}
   & \multicolumn{1}{r}{3.65 $\pm$ 1.20} \\
 & 4 
   & \multicolumn{1}{r}{0.92 $\pm$ 0.20}
   & \multicolumn{1}{r}{1.22 $\pm$ 0.30}
   & \multicolumn{1}{r}{5.81 $\pm$ 1.40} \\
 & 6 
   & \multicolumn{1}{r}{0.99 $\pm$ 0.20}
   & \multicolumn{1}{r}{1.32 $\pm$ 0.30}
   & \multicolumn{1}{r}{6.76 $\pm$ 1.50} \\
\midrule

\multirow{7}{*}{39}
 & 1  
   & \multicolumn{1}{r}{5.56 $\pm$ 0.60}
   & \multicolumn{1}{r}{5.74 $\pm$ 0.50}
   & \multicolumn{1}{r}{7.90 $\pm$ 1.10} \\
 & 2  
   & \multicolumn{1}{r}{5.66 $\pm$ 0.60}
   & \multicolumn{1}{r}{5.82 $\pm$ 0.70}
   & \multicolumn{1}{r}{10.60 $\pm$ 1.30} \\
 & 4  
   & \multicolumn{1}{r}{5.70 $\pm$ 0.60}
   & \multicolumn{1}{r}{6.23 $\pm$ 0.60}
   & \multicolumn{1}{r}{15.38 $\pm$ 1.30} \\
 & 8  
   & \multicolumn{1}{r}{5.80 $\pm$ 0.70}
   & \multicolumn{1}{r}{6.58 $\pm$ 0.70}
   & \multicolumn{1}{r}{21.20 $\pm$ 1.50} \\
 & 16 
   & \multicolumn{1}{r}{6.11 $\pm$ 0.70}
   & \multicolumn{1}{r}{7.07 $\pm$ 0.80}
   & \multicolumn{1}{r}{26.19 $\pm$ 1.80} \\
 & 32 
   & \multicolumn{1}{r}{6.14 $\pm$ 0.70}
   & \multicolumn{1}{r}{7.41 $\pm$ 0.80}
   & \multicolumn{1}{r}{29.70 $\pm$ 2.10} \\
 & 39 
   & \multicolumn{1}{r}{6.04 $\pm$ 0.60}
   & \multicolumn{1}{r}{7.31 $\pm$ 0.80}
   & \multicolumn{1}{r}{30.88 $\pm$ 2.30} \\
\midrule

\multirow{8}{*}{89}
 & 1  
   & \multicolumn{1}{r}{8.21 $\pm$ 1.00}
   & \multicolumn{1}{r}{8.51 $\pm$ 0.90}
   & \multicolumn{1}{r}{12.08 $\pm$ 1.40} \\
 & 2  
   & \multicolumn{1}{r}{8.33 $\pm$ 1.00}
   & \multicolumn{1}{r}{8.92 $\pm$ 1.00}
   & \multicolumn{1}{r}{15.70 $\pm$ 1.50} \\
 & 4  
   & \multicolumn{1}{r}{8.64 $\pm$ 1.00}
   & \multicolumn{1}{r}{9.61 $\pm$ 1.10}
   & \multicolumn{1}{r}{21.94 $\pm$ 1.80} \\
 & 8  
   & \multicolumn{1}{r}{9.09 $\pm$ 1.10}
   & \multicolumn{1}{r}{10.50 $\pm$ 1.20}
   & \multicolumn{1}{r}{29.82 $\pm$ 2.50} \\
 & 16 
   & \multicolumn{1}{r}{9.50 $\pm$ 1.20}
   & \multicolumn{1}{r}{11.53 $\pm$ 1.50}
   & \multicolumn{1}{r}{37.54 $\pm$ 2.60} \\
 & 32 
   & \multicolumn{1}{r}{9.75 $\pm$ 1.40}
   & \multicolumn{1}{r}{12.25 $\pm$ 1.70}
   & \multicolumn{1}{r}{42.49 $\pm$ 2.80} \\
 & 64 
   & \multicolumn{1}{r}{9.85 $\pm$ 1.30}
   & \multicolumn{1}{r}{12.75 $\pm$ 1.70}
   & \multicolumn{1}{r}{46.06 $\pm$ 2.60} \\
 & 89 
   & \multicolumn{1}{r}{9.96 $\pm$ 1.40}
   & \multicolumn{1}{r}{12.94 $\pm$ 1.90}
   & \multicolumn{1}{r}{46.30 $\pm$ 2.70} \\
\bottomrule
\end{tabular}
\end{table}

In Table~\ref{tab:my-table-percentage}, we report the relative errors with respect to the true minimum energies for Lennard-Jones systems of sizes \(N\in\{6,39,89\}\). Each row corresponds to a specific \((N,d)\) pair, where \(d\) is the number of coordinates updated per Metropolis-Hastings proposal, and each column indicates\( (\sigma_{\mathrm{total}}^{-2}\in\{200,100,10\}\). The entries list the mean relative error multiplied by 100, along with its standard deviation. 

We observe that, at fixed \(\sigma_{\mathrm{total}}^{-2}\), increasing \(d\) generally increases the relative error, indicating that spreading the same proposal variance over more coordinates tends to produce larger, less frequently accepted moves. Moreover, when \(\sigma_{\mathrm{total}}^{-2}\) is smaller, this effect becomes especially pronounced for large \(d\). Finally, smaller problems achieve consistently lower relative errors than larger systems, especially when both \(d\) and \(\sigma^2_{\mathrm{total}}\) are high. These results underscore the importance of tuning both the number of updated coordinates and the overall proposal variance to achieve accurate solutions in Lennard-Jones simulations.

\begin{figure}[htbp]
\centering
\includegraphics[width=0.9\textwidth]{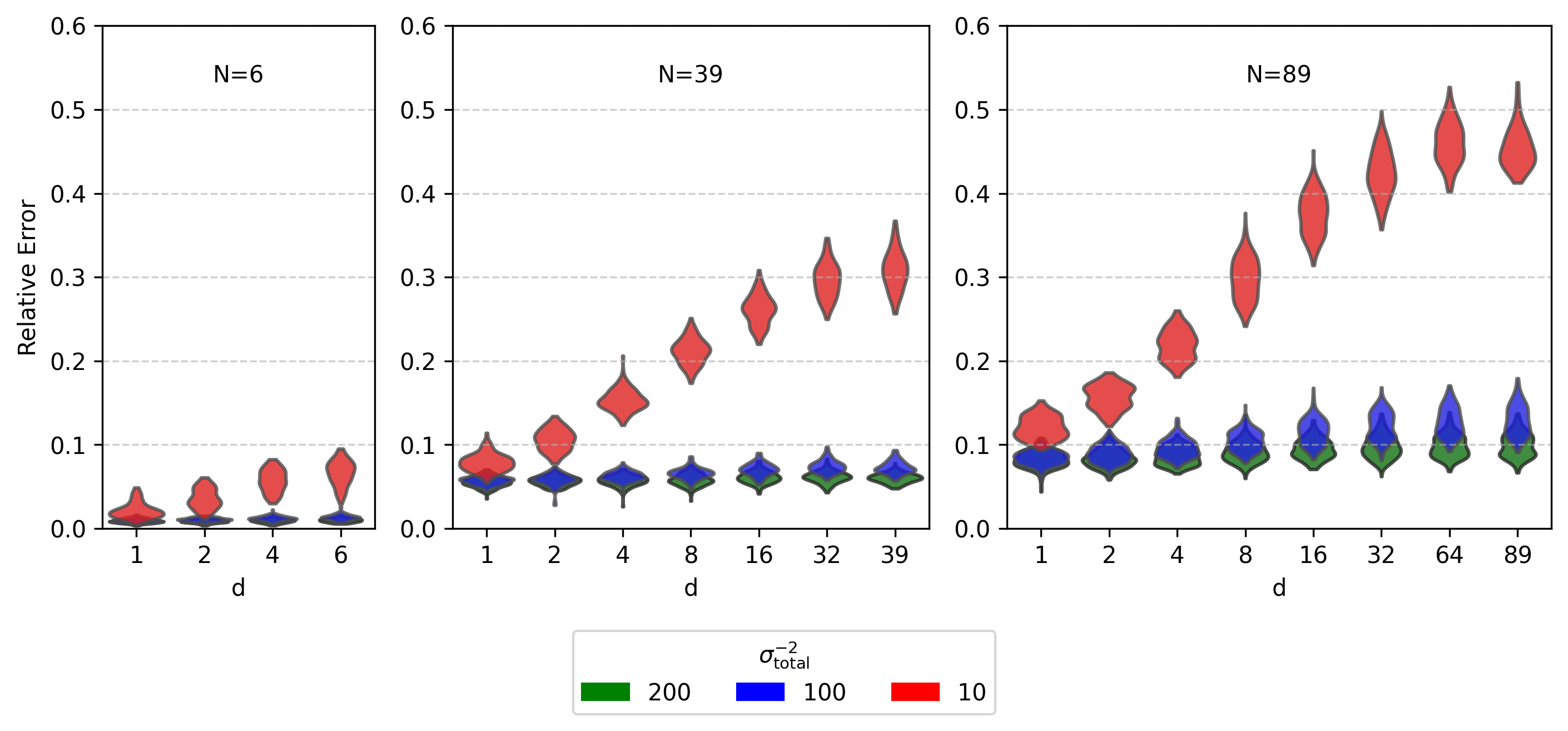}
\caption{Violin plots of the relative errors, 
\(
(\mathrm{Best} - \mathrm{BestMin}) / \lvert \mathrm{BestMin} \rvert
\),
for Lennard-Jones systems of sizes \(N=6,\,39,\,89\) at \(\text{Step}=100{,}000\). 
Each subplot has a custom width ratio of \(4{:}7{:}8\), and the horizontal axis 
shows \(d\) (the number of coordinates moved). 
Colors (red, green, blue) indicate three inverse-variance settings 
\(\sigma_{\mathrm{total}}^{-2}=200,\,100,\,10\). 
The vertical range is restricted to \([0,0.6]\), with labels 
\(\,N=6,\,39,\,89\) placed in the upper region of each subplot. 
A single figure-level legend at the bottom summarizes the three variance cases.}
\label{fig:OverlapViolin_DiffWidth}
\end{figure}

\subsubsection{Rosenbrock Problem}
\label{sub: experi:RP}

\noindent
\begin{figure}[t]
    \centering
    \includegraphics[width=0.99\textwidth]{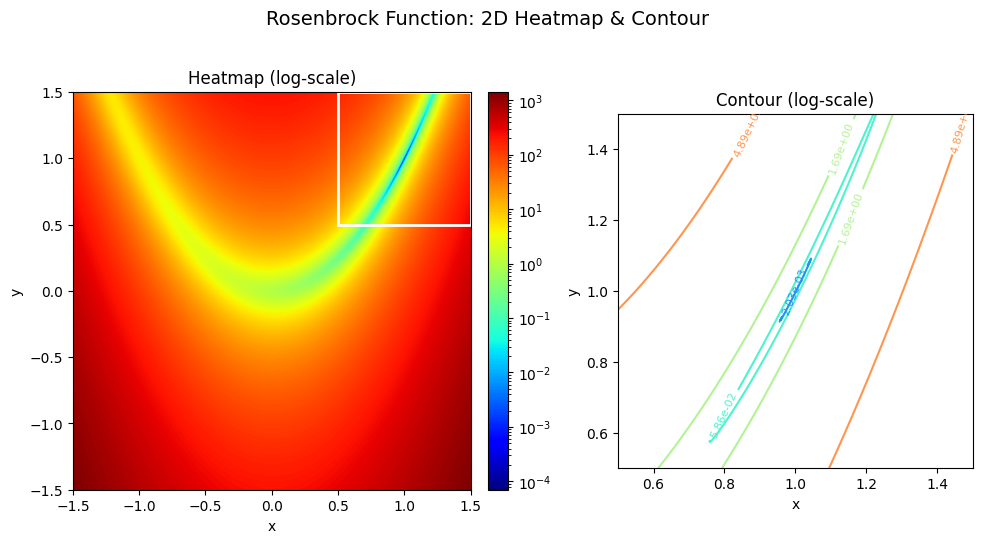} 
    \caption{Heatmap (left) and contour plot (right) for the 2D Rosenbrock function 
    $f(x,y) = (1 - x)^2 + 100\,\bigl(y - x^2\bigr)^2$ with logarithmic color scale. 
    The heatmap covers the domain $[-1.5,1.5]\times[-1.5,1.5]$, and a red rectangle 
    indicates the smaller region $[0.5,1.5]\times[0.5,1.5]$ for the contour plot.}
    \label{fig:rosenbrock-2d-heatmap-contour}
\end{figure}

We also tested the partial-coordinate Simulated Annealing (SA) schemes on the Rosenbrock function~\cite{rosenbrock1960automatic}, a classic benchmark in nonlinear optimization.  Recall that the Rosenbrock function in $N$ dimensions is given by
\[
  f(\mathbf{x}) 
  \;=\; 
  \sum_{i=1}^{N-1} \Bigl[100\,(x_{i+1} - x_i^2)^2 \;+\; (1 - x_i)^2 \Bigr],
  \quad
  \mathbf{x} = (x_1,\dots,x_N)\in\mathbb{R}^N,
\]
with a global minimum of $f(\mathbf{x})=0$ at $x_1 = x_2 = \cdots = x_N = 1$. This problem features a narrow ``banana-shaped'' valley, which is notoriously difficult for local methods.  


Table~\ref{tab:rosen_errors_final_sci} shows the final absolute errors across combinations of $\sigma_{\mathrm{total}}^{-2} \in \{6000,\,9000,\,12000\}$, partial‐update size $d$, and problem dimension $N = 30,\,72,\,200$. When $\sigma_{\text{total}}^{-2}$ is lower, errors remain small for $N=30$ and $72$, but can grow for $N=200$. As $\sigma_{\text{total}}^{-2}$ increases , the mean errors for $N=200$ can spike further, since more aggressive proposals are accepted less frequently later in the annealing schedule. Smaller update sizes generally show moderate errors in lower‐dimensional problems; however, when dimension and $\sigma_{\text{total}}^{-2}$ both grow, the errors can escalate. In contrast, moderate partial updates strike a balance between acceptance rate and step size, but with increasing variability for higher $N$. Overall, the data highlight that both the parameter $p$  and the inverse variance $\sigma_{\text{total}}^{-2}$ strongly influence the convergence pace and final accuracy on this classically ill‐conditioned problem.

\begin{table}[htbp]
\centering
\scriptsize  
\caption{Absolute errors for Rosen dataset across settings}
\label{tab:rosen_errors_final_sci}
\renewcommand{\arraystretch}{0.8}
\begin{tabular}{ccccc}
\toprule
\multicolumn{2}{c}{} & \multicolumn{3}{c}{\textbf{N}} \\
\cmidrule(lr){3-5}
\textbf{$1/\sigma^2_{\text{total}}$} & \textbf{$d$} & \textbf{30} & \textbf{72} & \textbf{200} \\
\midrule
\multirow{4}{*}{6000}
& 1 & 
$(5.64 \pm 8.01)\times 10^{-4}$ & 
$(8.62 \pm 10.14)\times 10^{-4}$ & 
$(9.58 \pm 17.99)\times 10^{-2}$ \\
& 2 & 
$(1.97 \pm 2.21)\times 10^{-3}$ & 
$(3.96 \pm 3.70)\times 10^{-3}$ & 
$(3.21 \pm 5.94)\times 10^{-1}$ \\
& 3 & 
$(5.35 \pm 5.48)\times 10^{-3}$ & 
$(9.04 \pm 5.78)\times 10^{-3}$ & 
$(5.52 \pm 10.60)\times 10^{-1}$ \\
& 6 & 
$(1.35 \pm 0.64)\times 10^{-2}$ & 
$(3.19 \pm 0.77)\times 10^{-2}$ & 
$(1.74 \pm 2.67)\times 10^{0}$ \\
\midrule
\multirow{3}{*}{9000}
& 1 & 
$(1.26 \pm 1.73)\times 10^{-3}$ & 
$(1.88 \pm 2.10)\times 10^{-3}$ & 
$(4.08 \pm 3.48)\times 10^{-1}$ \\
& 3 & 
$(8.42 \pm 7.21)\times 10^{-3}$ & 
$(2.10 \pm 1.30)\times 10^{-2}$ & 
$(1.49 \pm 2.31)\times 10^{0}$ \\
& 9 & 
$(3.72 \pm 1.16)\times 10^{-2}$ & 
$(1.00 \pm 0.22)\times 10^{-1}$ & 
$(2.53 \pm 4.92)\times 10^{0}$ \\
\midrule
\multirow{5}{*}{12000}
& 1 & 
$(1.69 \pm 2.08)\times 10^{-3}$ & 
$(1.89 \pm 2.22)\times 10^{-3}$ & 
$(7.36 \pm 4.55)\times 10^{-1}$ \\
& 2 & 
$(8.63 \pm 10.49)\times 10^{-3}$ & 
$(1.64 \pm 1.76)\times 10^{-2}$ & 
$(2.36 \pm 1.97)\times 10^{0}$ \\
& 3 & 
$(1.62 \pm 1.27)\times 10^{-2}$ & 
$(3.99 \pm 2.94)\times 10^{-2}$ & 
$(3.24 \pm 3.89)\times 10^{0}$ \\
& 4 & 
$(2.64 \pm 1.43)\times 10^{-2}$ & 
$(6.10 \pm 2.45)\times 10^{-2}$ & 
$(4.52 \pm 5.58)\times 10^{0}$ \\
& 6 & 
$(5.44 \pm 2.67)\times 10^{-2}$ & 
$(1.29 \pm 0.37)\times 10^{-1}$ & 
$(4.51 \pm 7.39)\times 10^{0}$ \\
\bottomrule
\end{tabular}
\end{table}

\begin{figure}[htbp]
    \centering
    \includegraphics[width=0.8\textwidth]{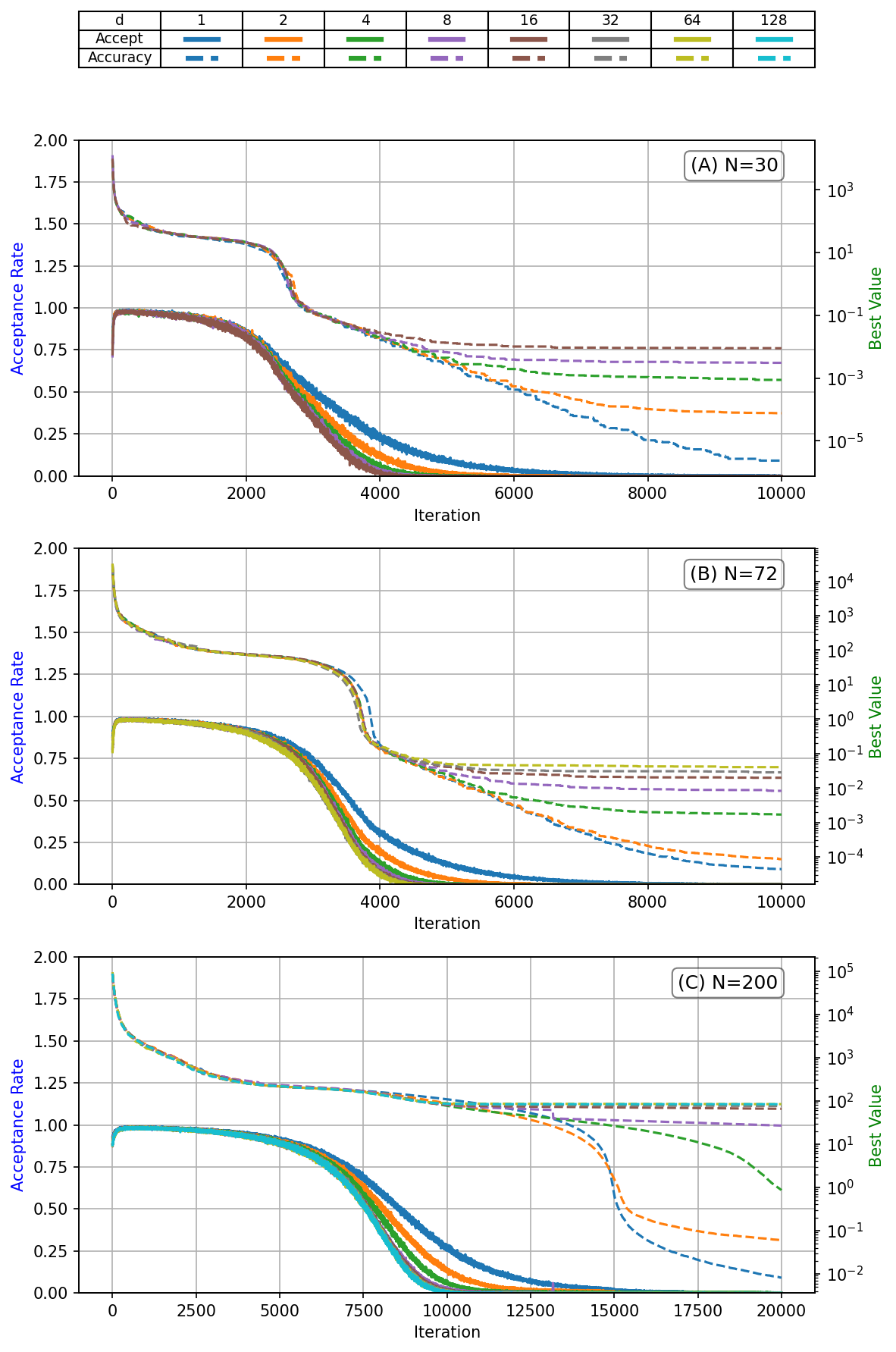}
    \caption{%
      Comparison across three problem sizes (N=30, 72, and 200).
      The legend at the top shows different \(d\) values for
      acceptance rate (solid lines) and best value (dashed lines).
    }
    \label{fig:one_large_figure}
\end{figure}

\subsubsection{Hyper-Elliptic proble}
\label{sub: experi:HEP}
\begin{figure}[t]
    \centering
    \includegraphics[width=0.99\textwidth]{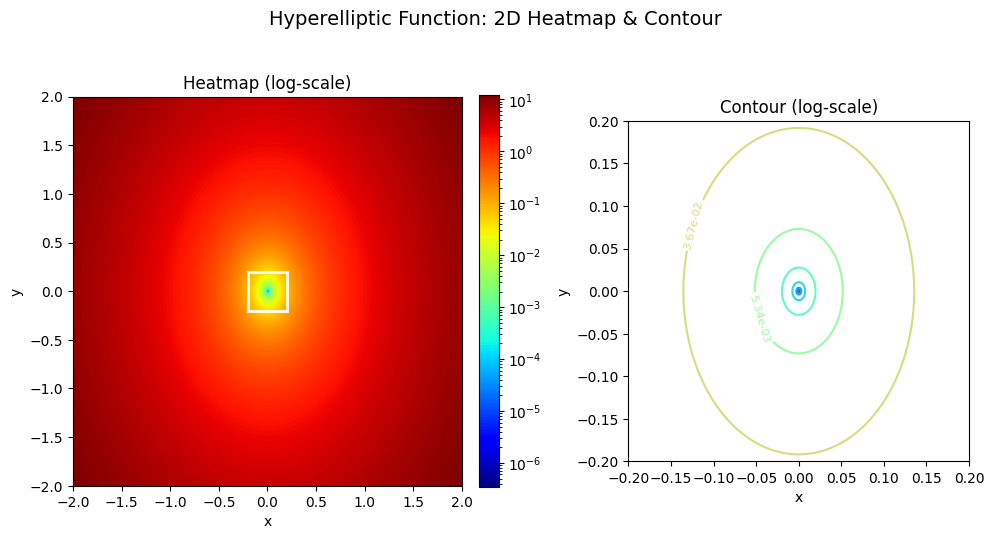} 
    \caption{Heatmap (left) and contour plot (right) for the hyperelliptic-like function 
    $f(x,y) = 2x^2 + y^2$ in 2D. 
    The heatmap spans the domain $[-2,2]\times[-2,2]$, and a red rectangle highlights 
    the smaller region $[-0.2,0.2]\times[-0.2,0.2]$ on which the contour plot is drawn.}
    \label{fig:hyperelliptic-2d-heatmap-contour}
\end{figure}

In this section, we examine a hyperelliptic-like function to illustrate how partial-coordinate simulated annealing (SA) performs in high-dimensional settings. Figure~\ref{fig:hyperelliptic-2d-heatmap-contour} offers a two-dimensional visualization of the function $f(x,y) = 2x^2 + y^2$, highlighting both its global structure (via a heatmap) and finer local contours (via a smaller zoomed-in region). By extending this idea to higher dimensions, we replicate the core challenge of navigating elongated valleys and “flat” directions that can slow convergence. We vary the partial‐update size $d$, the total proposal variance ($1/\sigma_{\text{total}}^2$), and the dimension $N$ to assess how robustly partial-coordinate SA handles this hyperelliptic‐like landscape. As shown in Table~\ref{tab:hyper_errors_final_sci}, tuning the combination of these parameters plays a crucial role in maintaining a balance between adequate exploration and sufficient acceptance rates, particularly when $N$ becomes large.

In our problem setting, the initial temperature was set to \(T_{0} = 2.0\). We then cooled the system to \(T = 0.2\) using a cooling factor of \(\alpha = 0.9999\). An equilibration phase of \(N_{\mathrm{equ}} = 100\) iterations was carried out before initiating a further cooling phase of \(N_{\mathrm{cooling}} = 200\) steps. Throughout these simulations, the system size was varied in \(N \in \{6, 39, 69\}\).

Table~\ref{tab:hyper_errors_final_sci} reports the final absolute errors for each combination of $1/\sigma_{\text{total}}^2 \in \{600,900,1200\}$, partial‐update size $d$, and problem dimension $N \in \{30,72,200\}$. As $\sigma_{\mathrm{total}}^{-2}$ increases, the proposals become bolder, which can help in escaping local traps but also risks diminishing acceptance rates at later stages. For lower dimensions, even relatively large $d$ values maintain small errors, whereas high $N$ exhibits more variation, especially when $d$ grows and $\sigma_{\mathrm{total}}^{-2}$ is large. In short, these results confirm that partial‐coordinate SA can effectively handle the hyperelliptic‐like constraints, though the interplay between partial‐update size and proposal variance must be tuned to balance acceptance rates with sufficiently rapid exploration.

\begin{figure}[htbp]
    \centering
    \includegraphics[width=0.8\textwidth]{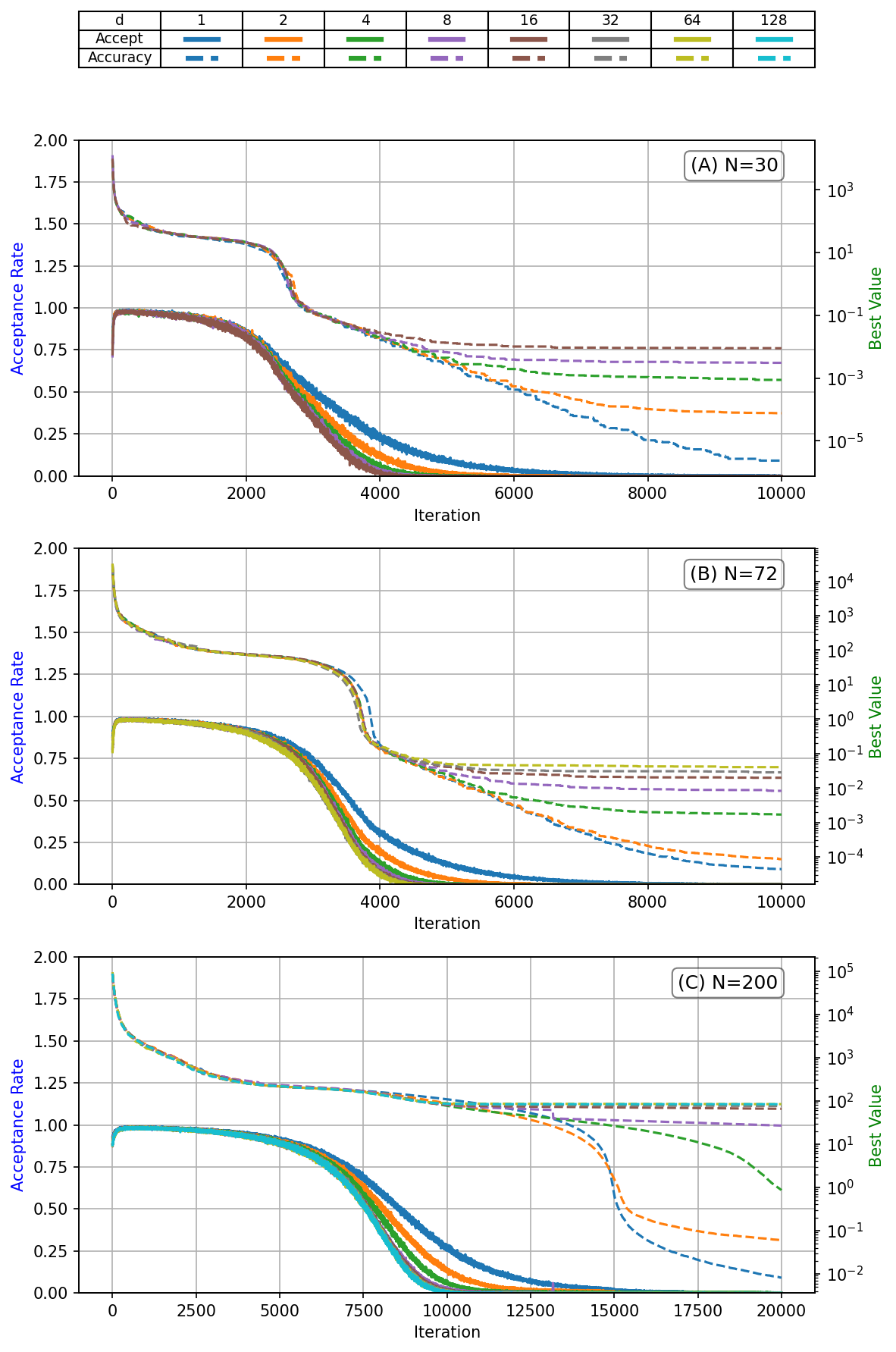}
    \caption{%
      Comparison across three problem sizes (N=30, 72, and 200).
      The legend at the top shows different \(d\) values for
      acceptance rate (solid lines) and best value (dashed lines).
    }
    \label{fig:one_large_figure}
\end{figure}
\begin{table}[htbp]
\centering
\footnotesize  
\caption{Absolute error for Hyper dataset across settings}
\label{tab:hyper_errors_final_sci}
\renewcommand{\arraystretch}{0.8}
\begin{tabular}{ccccc}
\toprule
\multicolumn{2}{c}{} & \multicolumn{3}{c}{$N$} \\
\cmidrule(lr){3-5}
\textbf{$1/\sigma^2_{\text{total}}$} & \textbf{$d$} & \textbf{30} & \textbf{72} & \textbf{200} \\
\midrule
\multirow{4}{*}{600}
& 1 
  & $(4.00 \pm 1.00)\times 10^{-6}$ 
  & $(3.20 \pm 0.60)\times 10^{-5}$ 
  & $(2.00 \pm 1.00)\times 10^{-6}$ \\
& 2 
  & $(1.42 \pm 0.32)\times 10^{-4}$ 
  & $(7.85 \pm 1.22)\times 10^{-4}$ 
  & $(7.20 \pm 0.76)\times 10^{-4}$ \\
& 3 
  & $(7.39 \pm 1.72)\times 10^{-4}$ 
  & $(3.19 \pm 0.46)\times 10^{-3}$ 
  & $(4.52 \pm 0.43)\times 10^{-3}$ \\
& 6 
  & $(3.20 \pm 0.51)\times 10^{-3}$ 
  & $(1.06 \pm 0.13)\times 10^{-2}$ 
  & $(2.12 \pm 0.15)\times 10^{-2}$ \\
\midrule
\multirow{3}{*}{900}
& 1 
  & $(6.00 \pm 2.00)\times 10^{-6}$ 
  & $(4.60 \pm 0.80)\times 10^{-5}$ 
  & $(5.00 \pm 1.00)\times 10^{-6}$ \\
& 3 
  & $(1.55 \pm 0.32)\times 10^{-3}$ 
  & $(6.62 \pm 0.84)\times 10^{-3}$ 
  & $(1.00 \pm 0.09)\times 10^{-2}$ \\
& 9 
  & $(9.77 \pm 1.49)\times 10^{-3}$ 
  & $(2.84 \pm 0.29)\times 10^{-2}$ 
  & $(6.20 \pm 0.39)\times 10^{-2}$ \\
\midrule
\multirow{5}{*}{1200}
& 1 
  & $(8.00 \pm 2.00)\times 10^{-6}$ 
  & $(5.90 \pm 1.20)\times 10^{-5}$ 
  & $(8.00 \pm 1.00)\times 10^{-6}$ \\
& 2 
  & $(4.50 \pm 1.04)\times 10^{-4}$ 
  & $(2.48 \pm 0.43)\times 10^{-3}$ 
  & $(2.69 \pm 0.28)\times 10^{-3}$ \\
& 3 
  & $(2.67 \pm 0.56)\times 10^{-3}$ 
  & $(1.14 \pm 0.17)\times 10^{-2}$ 
  & $(1.75 \pm 0.15)\times 10^{-2}$ \\
& 4 
  & $(5.97 \pm 1.15)\times 10^{-3}$ 
  & $(2.18 \pm 0.30)\times 10^{-2}$ 
  & $(4.09 \pm 0.33)\times 10^{-2}$ \\
& 6 
  & $(1.23 \pm 0.20)\times 10^{-2}$ 
  & $(4.00 \pm 0.50)\times 10^{-2}$ 
  & $(8.24 \pm 0.56)\times 10^{-2}$ \\
\bottomrule
\end{tabular}
\end{table}

\bigskip

\noindent


\section{Conclusions}
\label{sec:conclusions}
We have investigated the role of partial-coordinate moves in Metropolis--Hastings-based Simulated Annealing and demonstrated that selectively updating smaller subsets of coordinates can significantly improve performance in high-dimensional settings. Our theoretical analysis showed that partial-coordinate proposals help sustain a nontrivial acceptance rate when scaling the overall proposal variance, thereby maintaining larger effective step sizes and enhancing the chain’s mixing properties. In particular, we quantified how the acceptance probability and autocorrelation depend on the dimensionality of the proposal and derived approximate expressions for key MCMC metrics, including acceptance probabilities and third-order corrections via Edgeworth expansions. 

Extensive numerical experiments further confirmed that adjusting the number of coordinates \(d\) and distributing the total proposal variance judiciously can yield substantial gains in optimization accuracy. On Lennard-Jones clusters and high-dimensional benchmark problems such as Rosenbrock functions, partial-coordinate SA outperformed full-coordinate updates by achieving lower energies and faster convergence to near-optimal solutions, especially for larger problem dimensions. 

These results emphasize the importance of designing proposal mechanisms that carefully balance acceptance probability and exploration in high-dimensional Simulated Annealing. While we focused on Gaussian increments and isotropic proposals, the same principles can be extended to non-Gaussian or adaptive strategies. Future research directions include refining the subset selection procedure (for instance, using adaptive or problem-specific heuristics), incorporating multiple scales of step sizes within partial-coordinate updates, and applying these methods to more complex molecular or machine-learning models. Our findings suggest that partial-coordinate move strategies offer a flexible and efficient avenue for tackling large-scale optimization challenges via Simulated Annealing.



\end{document}